\documentclass[11pt,amsfonts]{article}
\usepackage{graphicx}
\usepackage{latexsym}
\usepackage{amssymb}
\usepackage{amsmath}
\usepackage{enumerate}
\usepackage{color}
\usepackage{layout}
\usepackage{dsfont}
\usepackage{eufrak}
\newtheorem{prop}{Proposition}
\newtheorem{lemma}{Lemma}

\newtheorem{theorem}{Theorem}
\newtheorem{remark}{Remark}

\def\real{{\mathord{{\rm I\kern-2.8pt R}}}}        
\def\inte{{\mathord{{\rm I\kern-2.8pt N}}}}

\def\sZZ{{\rm Z\kern-2.8ptem{}Z}}

\def\z{{\mathchoice
  {\sZZ}
  {\sZZ}
  {\rm Z\kern-0.30em{}Z}
  {\rm Z\kern-0.25em{}Z} }}
\def\sQQ{{\kern 0.27em \vrule height1.45ex width0.03em depth0em
          \kern-0.30em \rm Q}}
\def\qu{{\mathchoice
    {\sQQ}
    {\sQQ}
  {\kern 0.225em \vrule height1.05ex width0.025em depth0em \kern-0.25em \rm Q}
  {\kern 0.180em \vrule height0.78ex width0.020em depth0em \kern-0.20em \rm Q}
        }}
\def\sCC{{\kern 0.27em \vrule height1.45ex width0.03em depth0em
          \kern-0.30em \rm C}}
\def\complex{{\mathchoice
    {\sCC}
    {\sCC}
  {\kern 0.225em \vrule height1.05ex width0.025em depth0em \kern-0.25em \rm C}
  {\kern 0.180em \vrule height0.78ex width0.020em depth0em \kern-0.20em \rm C}
        }}

\newcommand{\ba}{\begin{array}}
\newcommand{\ea}{\end{array}}
\newcommand{\be}{\begin{equation}}
\newcommand{\ee}{\end{equation}}
\newcommand{\bea}{\begin{eqnarray}}
\newcommand{\eea}{\end{eqnarray}}
\newcommand{\beaa}{\begin{eqnarray*}}
\newcommand{\eeaa}{\end{eqnarray*}}

\def\z{\zeta}

\font\tenmath=msbm10 \font\sevenmath=msbm7 \font\fivemath=msbm5
\newfam\mathfam \textfont\mathfam=\tenmath
\scriptfont\mathfam=\sevenmath \scriptscriptfont\mathfam=\fivemath

\def \={{\buildrel {\rm (law)} \over =}}

\def\qed{ \hfill \vrule width.25cm height.25cm depth0cm\smallskip}

\newcommand{\basa}{\begin{assumption}}
\newcommand{\easa}{\end{assumption}}

\newcommand{\bas}{\begin{assum}}
\newcommand{\eas}{\end{assum}}

\newcommand{\ignore}[1]{}
\begin{document}

\renewcommand{\thefootnote}{\fnsymbol{footnote}}

\renewcommand{\thefootnote}{\fnsymbol{footnote}}

\title{ Existence and smoothness of the density for the  stochastic continuity equation  }
\author{David A.C. Mollinedo \footnote{Universidade Tecnol\'{o}gica Federal do Parana, Brazil. E-mail:  {\sl davida@utfpr.edu.br}
}
 $\quad$ Christian Olivera\footnote{Departamento de Matem\'{a}tica, Universidade Estadual de Campinas, Brazil.
E-mail:  {\sl  colivera@ime.unicamp.br}.  }
 $\quad$ Ciprian A.  Tudor \footnote{      Universit\'e de Lille 1. UFR Math\'ematiques.  E-mail:  {\sl Ciprian.Tudor@math.univ-lille1.fr}}
}
\maketitle
\begin{abstract}
We consider the stochastic  continuity equation driven by   Brownian motion. We
use the techniques of the Malliavin calculus  to show that the law of the solution has
a density with respect to the Lebesgue measure. We also
prove that the density  is H\"older continuous and satisfies some Gaussian-type estimates.
\end{abstract}

 \medskip


 \smallskip


\section{Introduction}

 In this paper we consider  the stochastic continuity equation given by

\begin{equation}\label{continuity}
 \left \{
\begin{aligned}
    &\partial_t u(t, x) + div  \left( ( b(t,x) + \frac{d B_{t}}{dt}) \cdot  u(t, x)  \right)= 0, \hskip0.3cm t\in [0,T], x\in \mathbb{R}^{d}
    \\[5pt]
    &u(0, x)=  u_{0} (x) \mbox{ for } x \in \mathbb{R}^{d} .
\end{aligned}
\right .
\end{equation}
where $(B_{t})_{t\in [0,T]}$ is a   Brownian motion (Bm) in $\mathbb{R}^d$  and the stochastic integration is understood in the Stratonovich sense. This equation constitutes a well-known model for several physical phenomena arising, among others,  in fluid dynamics and conservation of laws (see e.g. the monographs \cite{Li1}, \cite{Li2} or \cite{Daf}).
The stochastic continuity  equation  with random perturbation given by a  standard Brownian noise has been first studied in the celebrated works by Kunita \cite{Ku}, \cite{Ku3}. More recent works on these  topics are, among others,
\cite{Beck}, \cite{Fre1}, \cite{Gues} \cite{Moli} and \cite{NO}.

Our purpose is to analyze the density of the solution to the continuity equation in dimension $d=1$ by using Malliavin calculus. It is widely recognized today that this theory constitutes a powerful tool in order to study the densities of random variables in general, and of solutions to stochastic (partial) differential equations in particular. By using several criteria in terms of  the Malliavin derivative of the solution, we are able to show that the unique solution to (\ref{continuity}) admits a density which is H\"older continuous and satisfies some Gaussian estimates. These criteria, proved in \cite{Bally-Caramellino}, \cite{BH},  and \cite{NV} are recalled in the next section of our work.

In order to control the Malliavin derivative of the solution to (\ref{continuity}), we will use the representation theorem of the solution. It is well-known from \cite{Ku}, \cite{Ku3}  that one can associate to (\ref{continuity}) a so-called equation of characteristics whose solution generates a $C^{1}-$ flow of diffeomorphism. Then the  solution to the continuity equation  can be expressed as the initial condition applied to the inverse flow multiplied by the Jacobian of the inverse flow. Therefore, via a  careful analysis of the Malliavin derivatives of the  inverse flow, and by assuming suitable properties on the initial condition  $u_{0} $ and on the drift $b$  of the equation, we are able to control the Malliavin derivatives  of the solution and to prove several properties of the density of the law of the solution.

We organize our paper as follows. In Section 2, we present some preliminaries
on Malliavin calculus and the method of characteristics. In Section  3, we study the properties of the Malliavin derivative of the solution to (\ref{continuity}) while in section 4, we state and prove our main  results on regularity in law.

\section{Preliminaries}
In this preliminary section we present the basic elements from Malliavin calculus that we will need in this work.  In the second part we present three criteria based on Malliavin calculus that we will  use in order to prove the absolute continuity of the law of the solution with respect to the Lebesque measure and various properties of the density. The last paragraph contains some known facts on the continuity equation and on  its strong solution. In particular, we give the representation of the solution in terms of the initial condition and of the inverse flow which constitutes a key result for our approach.

\subsection{Notions of Malliavin Calculus}

	This subsection is devoted to present the basics tools from Malliavin calculus that will be used in the paper. We refer to 	\cite{N} for a complete exposition.  Consider ${\mathcal{H}}$ a real separable Hilbert space endowed with the scalar product $\langle \cdot, \cdot \rangle _{\mathcal{H}}$ and $(B (\varphi), \varphi\in{\mathcal{H}})$ an isonormal Gaussian process \index{Gaussian process} on a probability space $(\Omega, {\cal{A}}, \mathbb{P})$, that is, a centred Gaussian family of random variables such that $\mathbf{E}\left( B(\varphi) B(\psi) \right)  = \langle\varphi, \psi\rangle_{{\mathcal{H}}}$.
	
\medskip

We denote by $D$  the Malliavin  derivative operator that acts on smooth functions of the form $F=g(B(\varphi _{1}), \ldots , B(\varphi_{n}))$ ($g$ is a smooth function with compact support and $\varphi_{i} \in {{\cal{H}}}, i=1,...,n$)
\begin{equation*}
DF=\sum_{i=1}^{n}\frac{\partial g}{\partial x_{i}}(B(\varphi _{1}), \ldots , B(\varphi_{n}))\varphi_{i}.
\end{equation*}

It can be checked that the operator $D$ is closable from $\mathcal{S}$ (the space of smooth functionals as above) into $ L^{2}(\Omega; \mathcal{H})$ and it can be extended to the space $\mathbb{D} ^{1,p}$ which is the closure of $\mathcal{S}$ with respect to the norm
\begin{equation*}
\Vert F\Vert _{1,p} ^{p} = \mathbf{E} F ^{p} + \mathbf{E} \Vert DF\Vert _{\mathcal{H}} ^{p}.
\end{equation*}

We can analogously  define the $k$th iterated Malliavin derivative. We will denote by $\mathbb{D}^{k,p}$ with  $k\geq 2$ the completition of the set of smooth random variables with respect to the norm
\begin{equation*}
\Vert F\Vert ^{p} _{k,p} =\mathbf{E} F ^{p} + \sum_{j=1}^{k}\mathbf{E} \Vert DF\Vert _{\mathcal{H}^{\otimes j}} ^{p}
\end{equation*}

 We denote by  $ \mathbb{D} ^{k, \infty}:= \cap _{p\geq 1 } \mathbb{D} ^{k,p}$ for every $k\geq 1$. In this paper, $\mathcal{H}$ will be the standard Hilbert space $L^ {2}([0,T])$.

\medskip

We will use  the chain rule for the Malliavin derivative (see Proposition 1.2.4 in \cite{N}) which says that if $\varphi: \mathbb{R}\to \mathbb{R}$ is a differentiable function with bounded derivative and $F\in \mathbb{D} ^ {1,2}$, then  $\varphi (F) \in \mathbb{D} ^ {1,2}$ and
\begin{equation}
\label{chain}
D\varphi(F)= \varphi ' (F) DF.
\end{equation}
We also recall the rule to  differentiate a product of random variables: if $F,G \in \mathbb{D} ^ {1,2}$  such that $FG\in \mathbb{D} ^ {1,2}$, then
\begin{equation}
\label{prod}
D(FG)= FDG +GDF.
\end{equation}

\subsection{Criteria for densities of random variables}\label{three}

The Malliavin calculus is widely recognized theory as a powerful theory that allows to analyze the existence and  smoothness of the densities of certain random variables. We list below a collection of criteria  from Malliavin calculus that we will emply in our work.  The first is the well-known Bouleau-Hirsch criterium that says that the strict positivity of the norm of the  Malliavin derivative implies the absolute continuity of the law. The second criterium (from \cite{Bally-Caramellino}) gives sufficient conditions for the H\"older continuity of the density while the third criterium (given in \cite{NV}) allows to prove upper and lower Gaussian estimates for the density of a random variable.

\medskip

A notable  result is the Bouleau-Hirsch theorem, which provides a criteria in terms of Malliavin derivatives for a random variable to have a density (see \cite{BH} or  \cite{N}).

\begin{theorem}\label{Baleu}
Let $F$ be a random variable of the space $\mathbb{D}^{1,2}$ and suppose that $\|DF \|_{L^{2}([0,T])} > 0$ a.s. Then the law of $F$ is absolutely continuous with respect to the Lebesgue measure on $\mathbb{R}$.
\end{theorem}

In \cite{Bally-Caramellino},  Proposition 23,  the authors  have presented sufficient conditions for a random variable to have a H\"older continuous density.

\begin{prop} \label{Prob_B-C}
	Let $F \in \cap_{p\in \mathbb{N}} \mathbb{D}^{2,p}$ and assume that $\Vert DF\Vert _{\mathcal{H}}^{-2} \in \cap_{p\in \mathbb{N}} L^{p}(\Omega)$. Then
the law of $F$ is absolutely continuous with respect to the Lebesgue measure on $\mathbb{R}$ and its density $\rho_{F}$ is H\"older continuous of any exponent $0< q < 1$. Moreover,  there exists some universal constant $C>0$ and $p_{i}>0$, $i = 1, \ldots , 5$, such that
$$
\rho_{F}(z) \leq C \mathbf{E}\left( \Vert DF\Vert _{\mathcal{H}}^{- 2 p_{1}} \right)^{p_{2}} \Vert F \Vert_{2, p_{4}}^{p_{3}} \left( \mathbb{P}(| F - z | \leq 2 ) \right)^{1/p_{5}}.
$$
In particular, $\lim _{z \rightarrow \infty} |z|^{p} \rho_{F}(z) = 0$ for every $p \in \mathbb{N}$.
\end{prop}

The main tool in order to obtain the Gaussian estimates for the density of the solution to the continuity equation is the following result given in \cite{NV}, Corollary 3.5.

\begin{prop}\label{6a-4}

If $ F \in \mathbb{D} ^{1,2}$, let
$$g_{F}(F)= \int_{0} ^{\infty} d\theta e^{-\theta} \mathbf{E}  \left[ \mathbf{E}  '\left( \langle DF, \widetilde{DF} \rangle _{\mathcal{H}} | F \right) \right]$$
where for any random variable $X$, we denote
\begin{equation}
\label{tilde}
\tilde{X}(\omega , \omega ') = X ( e ^{-\theta }w + \sqrt{1-e^{-2\theta} }\omega ').
\end{equation}
Here  $\tilde{X}$ is  defined on a product probability space
 $\left( \Omega \times \Omega ', {\cal{F}} \otimes {\cal{F}}, P\times P'\right)$ and ${\mathbf E}'$
 denotes the expectation with respect to the probability measure $P'$.
If there exist two constants $\gamma _{min}>0$ and $\gamma _{max}>0$ such that  almost surely
$$0\leq \gamma_{ min}^{2} \leq g_{F}(F) \leq \gamma _{max} ^{2}$$
then $F$ admits a density $\rho$. Moreover, for every $z\in \mathbb{R}$,
\begin{equation*}
\frac{ \mathbf{E} \vert F -\mathbf{E} F \vert }{2\gamma ^{2}_{max} }e ^{ -\frac{(z-\mathbf{E} F) ^{2}}{2\gamma^{2} _{min} } }\leq \rho(z) \leq
\frac{ \mathbf{E} \vert F -\mathbf{E} F \vert }{2\gamma^{2} _{min} }e ^{ -\frac{(z-\mathbf{E} F) ^{2}}{2\gamma ^{2} _{max} }}.
\end{equation*}

\end{prop}

\subsection{Representation of the solution}

We will start by recalling some known facts on the solution to the  stochastic continuity  equation driven by a standard Wiener process in $\mathbb{R}^{d}$.

\medskip

The equation (\ref{continuity}) is interpreted in the strong sense, as the solution to the following stochastic
integral equation
\begin{equation}
u(t,x)=u_{0}(x)-\int_{0}^{t}  div (b(s,x) u(s,x) ) \ ds -\sum_{i=0}^{d}\int_{0}%
^{t}   \partial_{x_i}  u(s,x)\circ dB_{s}^{i} \label{transintegral}
\end{equation}
for $t\in [0,T]$ and $x\in \mathbb{R} ^{d}$.

\medskip

The solution to (\ref{continuity})  is related with the so-called equation of characteristics. That is, for $0\leq s\leq t$ and $x\in\mathbb{R}^{d}$, consider the following
stochastic differential equation in $\mathbb{R}^{d}$

\begin{equation}
\label{11}X_{s,t}(x)= x + \int_{s}^{t} b(r, X_{s,r}(x)) \ dr + B_{t}-B_{s},
\end{equation}
and denote by $X_{t}(x): = X_{0,t}(x), t\in [0,T], x\in \mathbb{R} ^{d}$.

\medskip

For $m \in \mathbb{N}$ and $0< \alpha < 1$, let us assume the following hypothesis on $b$:
\begin{equation}
\label{REGULCLASS}
    b\in L^{1}((0,T); C_{b}^{m,\alpha}(\mathbb{R}^{d}))
\end{equation}
where $C^{m,\alpha}(\mathbb{R}^{d})$ denotes the class of functions of class $C ^{m}$  on $\mathbb{R}^{d}$ such that the last derivative is H\"older continuous of order $\alpha$. By considering the condition (\ref{REGULCLASS}) it is well known that $X_{s,t}(x)$ is a
stochastic flow of $C^{m}$-diffeomorphism (see for example \cite{Chow} and
\cite{Ku}). Moreover, the inverse flow
 $$Y_{s,t}(x):=X_{s,t}^{-1}(x)$$
 satisfies the
following backward stochastic differential equation%

\begin{equation}
\label{back}Y_{s,t}(x)= x - \int_{s}^{t} b(r, Y_{r,t}(x)) \ dr - (B_{t}-B_{s}),
\end{equation}
for  every $0\leq s\leq t\leq T$. We will denote $Y_{0,t}(x):= Y_{t}(x)$ for every $t\in [0, T], x\in \mathbb{R}$.

\medskip

We have the following representation of the solution to the  stochastic continuity  equation in terms of the inital data and of the inverse flow  (\ref{back}). We refer to e.g. \cite{Ku}  or \cite{Chow}, Section 3 for the proof.  By $J$ we denote the Jacobian.

\begin{lemma}
\label{lemaexis}  Assume (\ref{REGULCLASS}) for $m\geq3, \delta >0$ and let $u_{0}\in C^{m,\delta}(\mathbb{R}^{d})$. Then the Cauchy problem (\ref{transintegral}) has a
unique solution $(u(t,x))_{t\in [0,T], x\in \mathbb{R}^{d}} $ which can be represented as

\begin{equation}\label{rep}
u(t,x)= u_{0}(Y_{t}(x)) JY_{t}(x),  \hskip0.3cm t \in [0,T], x\in \mathbb{R} ^{d}.
\end{equation}
\end{lemma}

Our approach to prove the regularity of the solution in the sense of Malliavin calculus is based on the formula (\ref{rep}) by assuming that the initial condition is smooth. A similar representation to (\ref{rep}) exists for the transport equation (see e.g.  \cite{Ku}) and it has been used in \cite{OT} to obtain to absolute continuity of the law of the solution to the transport equation.

\section{Properties of the inverse flow and Malliavin differentiability of the solution}

Taking into account the formula (\ref{rep}), in order to control the Malliavin derivative of $u(t,x)$, one needs to analyse the inverse flow.   We will need to control   its Malliavin derivative, its Jacobian and the Malliavin derivative of its Jacobian. We will restrict, from now on, to the case $d=1$.

\medskip

Consider  $  \left( Y_{s,t}(x)\right)_{ 0\leq s\leq t\leq T, x \in \mathbb{R}}$  the inverse flow $ Y _{s,t}(x)= X_{s,t}^{-1}(x)$, $0\leq s\leq t\leq T, x\in \mathbb{R}$ with $X_{s,t}(x)$ from (\ref{11}).
Recall that satisfies the
following backward stochastic differential equation (\ref{back}).

\medskip

We will keep the multidimensional   notation $JY _{s,t}(x)= \frac{d}{dx} Y_{s,t}(x)$ although we work now in dimension $d=1$. We also denote by $b'(t,x)= \frac{d}{dx } b(t,x) $ , $b''(t,x)= \frac{d^{2}}{dx^{2}  } b(t,x) $ and  $ b^{(n) }(t,x) = \frac{ d ^{n}}{dx^{n} }b(t,x)$.
By $\Vert f\Vert _{\infty} $ we denote the infinity norm of the function $f$.

\vskip0.2cm

Concerning the Malliavin derivative of the inverse flow, we have the following estimate.

\begin{lemma} \label{l1}Assume $ b\in L ^{\infty} \left( (0,T), C^{1} _{b}(\mathbb{R})\right)$.

If $  \left( Y_{s,t}(x)\right)_{ 0\leq s\leq t\leq T, x \in \mathbb{R}}$ denotes the inverse flow (\ref{back}), then for every  $0\leq s\leq t\leq T$ and $x\in \mathbb{R}$ the random variable $ Y_{s,t}(x)$ is Malliavin differentiable  and we have,  for every $\alpha \in (0, T]$
\begin{equation}\label{5a-4}
D_{\alpha} Y _{s,t} (x)=- 1_{[s,t]} (\alpha )e^{- \int_{s} ^{\alpha} b'(r, Y_{r,t}) dr}.
\end{equation}
Moreover, there exists a constant $ C_{1}>0$ (not depending on $\omega$) such that
\begin{equation}\label{5o-1}
\sup_{\alpha \in (0, T), 0\leq s\leq t\leq T} \sup_{x \in\mathbb{R}}\vert D_{\alpha} Y _{s,t} (x)\vert \leq C_{1}.
\end{equation}
\end{lemma}
{\bf Proof: } The Malliavin differentiability of $ Y_{s,t}(x)$ and the formula (\ref{5a-4}) have been proven in Proposition 2 in  \cite{OT}. We also have, for every $\alpha, s, t, x$,
\begin{equation*}
\vert D_{\alpha} Y _{s,t} (x)\vert \leq e^{- \int_{s} ^{\alpha} b'(r, Y_{r,t}) dr} \leq e ^{ T\Vert b'\Vert _{\infty} }:= C_1
\end{equation*}
by using the fact that $b'$ is uniformly  bounded. This  implies (\ref{5o-1}).  \qed

\medskip

Now, we regard the Jacobian of $Y$.

\begin{lemma}\label{l2}Assume $ b\in L ^{\infty} \left( (0,T), C^{1} _{b}(\mathbb{R})\right)$.
Consider  $  \left( Y_{s,t}(x)\right)_{ 0\leq s\leq t\leq T, x \in \mathbb{R}}$ given by (\ref{back}).  Then,  for every  $0\leq s\leq t\leq T$ and $x\in \mathbb{R}$,

\begin{equation}\label{jy}
JY_{s,t}(x)=e ^{-\int_{s}^{t} b'(v, Y_{v,t}(x) )dv}
\end{equation}
and there exists $C_1:=e^{ T\Vert b'\Vert _{\infty} }$ such that
\begin{equation}\label{5o-2}
\sup_{0\leq s\leq t\leq T} \sup_{x \in\mathbb{R}}\vert JY _{s,t} (x)\vert \leq C_{1}.
\end{equation}
\end{lemma}
{\bf Proof: } By differentiating with respect to $x$ in (\ref{back}), we get
\begin{equation*}
J Y_{s,t}(x)= 1- \int_{s}^{t} b'(r, Y_{r,t}(x)) JY _{r,t}(x) dr
\end{equation*}
and this implies (\ref{jy}). Since $ b\in L ^{\infty} \left( (0,T), C^{1} _{b}(\mathbb{R})\right)$, we have
\begin{equation*}
 \vert JY_{s,t}(x)\vert \leq e ^{ T\Vert b'\Vert _{\infty} }
\end{equation*}
for every $s,t,x$ and  (\ref{5o-2}) follows. \qed

\medskip

We have a similar result for the Malliavin derivative of the inverse flow. We need in addition to assume the existence and boundness of the second derivative of the drift $b$.

\begin{lemma}\label{l3}Assume $ b\in L ^{\infty} \left( (0,T), C^{2} _{b}(\mathbb{R})\right)$.
Then $JY_{s,t}(x)$ belongs to $\mathbb{D}^{1,p} $ for every  $0\leq s\leq t\leq T$ and $x\in \mathbb{R}$ and there exists $C_{2}>0$ with
\begin{equation}\label{5o-3}
\sup_{\alpha \in (0, T), 0\leq s\leq t\leq T} \sup_{x \in\mathbb{R}}\vert D_{\alpha} JY _{s,t} (x)\vert \leq C_{2}.
\end{equation}
\end{lemma}
{\bf Proof: } From (\ref{jy}), it is clear that $JY_{s,t}(x)$ belongs to $\mathbb{D} ^{1,p} $ for every  $0\leq s\leq t\leq T$ and $x\in \mathbb{R}$. By using the chain rule (\ref{chain}),
\begin{equation}\label{6o-4}
D_{\alpha} JY _{s,t}(x)= -e ^{-\int_{s}^{t} b'(v, Y_{v,t}(x) )dv} \int_{s}^{t} b''(v, Y_{v,t}(x) )D_{\alpha } Y _{v,t} (x) dv
\end{equation}
and by (\ref{5o-1}), we obtain
\begin{equation*}
\vert D_{\alpha} JY _{s,t} (x)\vert  \leq  C_{1}^2 T  \Vert b''\Vert_{\infty} :=C_{2}.
\end{equation*}
\qed

\medskip

From the three lemmas above and the representation (\ref{rep}), we immediately obtain the following result on the Malliavin differentiability of the solution to the continuity equation.

\begin{prop}\label{p2}
 Let $\left( u(t,x)\right) _{t\in [0, T], x\in \mathbb{R}^{d}} $ be given by (\ref{transintegral}).
\begin{enumerate}

\item Assume $u_{0} \in C_{b}^{1} (\mathbb{R}) $ and $ b\in L^{\infty} \left( (0,T), C^{2} _{b}(\mathbb{R})\right)$.  Then $u(t,x)$ belong to $\mathbb{D}^{1,p} $ for every $t,x$ and for every $p\geq 2$ and there exists  $C_{3}>0$ such that
 \begin{equation}\label{5o-4}
\sup_{\alpha \in (0, T), 0\leq t\leq T} \sup_{x \in\mathbb{R}}\vert D_{\alpha} u(t,x)\vert \leq C_{3}.
\end{equation}
\item Assume $u_{0} \in C_{b}^{2} (\mathbb{R}) $ and $ b\in L ^{\infty} \left( (0,T), C^{3}_{b}(\mathbb{R})\right)$. Then $u(t,x)$ belong to $\mathbb{D} ^{2,p} $ for every $t,x$ and for every $p\geq 2$ and there exists  $C_{4}>0$ such that
 \begin{equation}\label{5o-12}
\sup_{\alpha, \beta \in (0, T)} \sup_{0\leq t\leq T}\sup_{x \in\mathbb{R}}\vert D_{\beta } D_{\alpha} u(t,x)\vert \leq C_{4}.
\end{equation}

\end{enumerate}
\end{prop}
{\bf Proof: } 1. We start by proving the first point. From (\ref{rep}) and Lemmas \ref{l1}, \ref{l2}, \ref{l3} it is clear that the random variable $u(t,x)$ is differentiable in the Malliavin sense for all $t, x$. Using the product rule (\ref{prod})  and the chain rule (\ref{chain})  for the Malliavin derivative (recall that $Y_{0,t}= Y_{t}$),
\begin{eqnarray}
D_{\alpha} u(t,x)&=&D_{\alpha} u_{0}(Y_{t}(x))  JY_{t}(x)+ u_{0}(Y _{t}(x))D_{\alpha} JY_{t}(x)\nonumber \\
&=& u_{0} ' (Y_{t}(x))D_{\alpha} Y _{t}(x) JY_{t}(x)+u_{0}(Y _{t}(x))D_{\alpha} JY_{t}(x).\label{5o-5}
\end{eqnarray}
From the estimates (\ref{5o-1}), (\ref{5o-2}) and (\ref{5o-3}) and the assumption on $u_{0}$, we obtain
\begin{equation*}
\sup_{\alpha \in (0, T), 0\leq t\leq T} \sup_{x \in\mathbb{R}}\vert D_{\alpha} u(t,x)\vert \leq \Vert u_{0} '\Vert _{\infty}  C_{1}^2 + \Vert u_{0} \Vert _{\infty} C_{2}:=C_{3}.
\end{equation*}
The above relation also implies that $u(t,x)\in \mathbb{D}^{1,p} $ for every $p\geq 2$.

2. Concerning the second point, we note that obviously $  D_{\alpha} u(t,x) $ is also Malliavin differentiable and from (\ref{5o-5}),
\begin{eqnarray*}
D_{\beta} D _{\alpha} u(t,x) &=&u_{0} ' (Y_{t}(x)) D_{\alpha} Y_{t}(x) D_{\beta} JY_{t}(x) + u_{0} ' (Y_{t}(x))D_{\beta} D_{\alpha}  Y_{t}(x) JY_{t}(x) \\
&&+u_{0} ''  (Y_{t}(x)) D_{\beta} Y_{t}(x)D _{\alpha} Y_{t}(x) JY_{t}(x)+ u_{0} (Y_{t}(x))D_{\beta} D_{\alpha} JY_{t}(x) \\
&&+ u_{0} ' (Y_{t}(x))D_{\beta} Y_{t} (x) D_{\alpha} JY_{t}(x).
\end{eqnarray*}
We need to bound the second Malliavin derivative of $Y_{t}(x)$ and $JY_{t}(x)$. We have from (\ref{5a-4}) (remember that $Y_{0,t}= Y_{t}$)
\begin{equation}\label{6o-1}
D_{\beta} D _{\alpha} Y_{t}(x) =1_{[0,t]}(\alpha) e^{-\int_{0}^{t} b'(v, Y_{v,t}(x)) dv} \int_{0}^{t} b'' (v, Y_{v,t}(x) ) D_{\beta} Y_{v, t}(x) dv
\end{equation}
so
\begin{equation}\label{6o-5}
\left| D_{\beta} D_{\alpha} Y_{t}(x)\right|  \leq e^{T \Vert b' \Vert_{\infty}}T\Vert b''\Vert_{\infty} C_{1}=C_1^2 T\Vert b'' \Vert_{\infty} = C_2.
\end{equation}
Also, from (\ref{6o-4}), for $\alpha , \beta \in (0,T]$
\begin{eqnarray*}
&&D_{\beta} D_{\alpha} JY_{t}(x)\\
&&= -e ^{-\int_{0}^{t} b'(v, Y_{v,t}(x) )dv}\int_{0}^{t} \Big(  b ^ {(3)} (v, Y_{v,t}(x) )D_{\beta} Y_{v,t}(x) D_{\alpha}Y_{v,t}(x) + b'' (v, Y_{v,t} (x))D_{\beta} D _{\alpha} Y_{v,t}(x)\Big) dv \\
&&\,\,\,\,\,+ e^ {-\int_{0}^{t} b'(v, Y_{v,t}(x) )dv } \left( \int_{0}^{t} b'' (v, Y_{v,t}(x) ) D _{\beta} Y_{v,t}(x))dv\right) \left( \int_{0}^{t} b'' (v, Y_{v,t}(x) ) D _{\alpha} Y_{v,t}(x))dv\right).
\end{eqnarray*}
Considering the estimates (\ref{6o-5}) and (\ref{5o-1}), and the assumption on $b$, we obtain
\begin{eqnarray*}
\left| D_{\beta} D_{\alpha} JY_{t}(x)\right| & \leq & e^{T \Vert b' \Vert _{\infty}} \Big[ T C_1^2\|b^{(3)}\|_{\infty} + T C_2 \|b''\|_{\infty}+ T^2 C_1^2 \|b''\|_{\infty}^2\Big] \\
&\leq&  T C_1^3 \Big[ \|b^{(3)}\|_{\infty}+2 T \|b''\|^2_{\infty}\Big].
\end{eqnarray*}
Hence, from the last two estimates, together with (\ref{5o-1}), (\ref{5o-2}) and (\ref{5o-3}), imply the conclusion of point 2. \qed

\section{Density of the solution}
In this section, we prove the existence and various properties  of the density of the solution to (\ref{continuity}) via the three criteria presented in Section \ref{three}. Under a relatively strong control of the initial condition and of the drift, we will show that these criteria can be applied to (\ref{transintegral}).

\subsection{Existence of the density}

The first result gives the absolute continuity of the law of the solution when the second derivative of the drift is negative and $u_{0}, u_{0} '$ are bounded and positive.

\begin{prop}\label{existence}  We will assume that $ b\in L ^{\infty} \left( (0,T), C^{2} _{b}(\mathbb{R})\right)$ and

\begin{equation}
\label{cc1}
b''(t,x) \leq 0 \mbox{ for all }t\in [0,T], x\in \mathbb{R}.
\end{equation}
Assume $u_{0} \in C_{b} ^{1} (\mathbb{R}) $ satisfies
\begin{equation}
\label{cc2}
u_{0}(x) > 0 \mbox{ and } u_{0} '(x) >0 \mbox{ for every } x \in \mathbb{R}.
\end{equation}
 If $(u(t,x))_{t\in [0, T], x\in \mathbb{R}}$ is the  unique  solution to the stochastic continuity equation \eqref{continuity}, then, for every $t\in (0,T]$ and for every $x\in \mathbb{R}$,
the law of $u(t,x)$   is absolutely continuous with respect to the Lebesgue
measure  on $\mathbb{R}$.
\end{prop}
{\bf Proof: } We need to show that $\Vert Du(t,x)\Vert _{\mathcal{H}}>0$ almost surely where $\mathcal{H}= L ^{2} ([0,T])$. From (\ref{5o-5}), we have
\begin{eqnarray*}
\Vert Du(t,x)\Vert  ^{2} _{\mathcal{H}}&=&\left( u_{0} ' (Y_{t}(x))\right) ^{2} (JY_{t}(x)) ^{2} \Vert DY_{t}(x)\Vert_{\mathcal{H}} ^{2} + (u_{0} (Y_{t}(x)) ) ^{2}  \Vert DJY_{t}(x)\Vert_{\mathcal{H}} ^{2} \\
&&+ 2u_{0} ' (Y_{t}(x))JY_{t}(x)u_{0}(Y _{t}(x)) \langle DY_{t}(x), DJY_{t}(x) \rangle _{\mathcal{H}}\\
&:=& A_{1}+A_{2}+A_{3}.
\end{eqnarray*}
From (\ref{cc2}), by noticing that $(JY_{t}(x)) ^{2}$ and  $\Vert DY_{t}(x)\Vert _{\mathcal{H}} ^{2}$ are strictly positive due to (\ref{5a-4}) and  (\ref{jy}), we have $A_{1}>0$. Clearly $A_{2}\geq 0.$  Concerning the summand $A_{3}$, by the condition (\ref{cc1}), we have from (\ref{5a-4}) and (\ref{6o-4}), for $\alpha \in (0, T]$, that
\begin{eqnarray}\label{5o-6}
D_{\alpha} JY _{s,t}(x)&=&
 -e ^{-\int_{s}^{t} b'(v, Y_{v,t}(x) )dv} \int_{s}^{t} b''(v, Y_{v,t}(x) )D_{\alpha } Y_{v,t} (x) dv \nonumber\\
&=&e ^{-\int_{s}^{t} b'(v, Y_{v,t}(x) )dv} \int_{s}^{t} b''(v, Y_{v,t}(x) )1_{[v,t]} (\alpha )e^{- \int_{v} ^{\alpha} b'(r, Y_{r,t}) dr} dv \\
&\leq& 0\,. \nonumber
\end{eqnarray}
Since $D_{\alpha} Y_{t}(x)\leq 0$ for every $\alpha, t\in (0,T], x\in \mathbb{R} $  due to (\ref{5a-4}) we have $\langle DY_{t}(x), DJY_{t}(x) \rangle_{\mathcal{H}}\geq 0$.  By using the assumption (\ref{cc2}) we get $A_{3} \geq 0$. Therefore $\Vert Du(t,x)\Vert _{\mathcal{H}}>0$ and the conclusion follows from  Theorem \ref{Baleu}. \qed

\begin{remark}
For example, the function $u_{0} (x)= \frac{\pi}{2} + \arctan (x)+\delta $ with $\delta >0$ satisfies the assumption in Proposition \ref{existence} (and also in Lemma \ref{l6} and Propositions \ref{regularity},  \ref{p6} below).
\end{remark}

\subsection{H\"older continuity of the density}
We apply now Proposition \ref{Prob_B-C} in order to show that the density of $u(t,x)$ is H\"older continuous. To this end, we will also need to control the second derivative of the solution.  We start by proving that $\Vert Du(t,x) \Vert ^{2}_{\mathcal{H}}  $ is bounded below by a positive constant which does not depend on $\omega$.

\begin{lemma}\label{l6}
Assume  $b\in L ^{\infty} \left( (0,T), C^{2} _{b}(\mathbb{R})\right)$ such that
\begin{equation}
\label{cc11}
b''(t,x) \leq -C<0 \mbox{ for all }t\in [0,T], x\in \mathbb{R}.
\end{equation}
Assume $u_{0} \in C_{b} ^{1} (\mathbb{R}) $ satisfies
\begin{equation}
\label{cc22}
u_{0}(x) \geq C > 0 \mbox{ and } u_{0} '(x) > 0 \mbox{ for every } x \in \mathbb{R}.
\end{equation}
Then  for $t\in (0,T]$ and $x\in \mathbb{R}$, there exists a constant $C_5=C_5(t)>0$ (non depending on $\omega$ and $x$) such that
$$\Vert Du(t,x) \Vert ^{2}_{\mathcal{H}}  \geq C_5(t)$$
\end{lemma}
{\bf Proof: } Notice that from (\ref{5a-4}), for $t, \alpha \in (0, T], x\in \mathbb{R}$, and by the assumption \eqref{cc11} we have
\begin{eqnarray*}
&\Vert DJY_{t}(x)\Vert_{\mathcal{H}} ^{2}&=e^{-2 \int_{0}^{t} b'(v, Y_{v,t}(x) )dv}\Big\langle  \int_0^t b''(Y_{v,t}(x))D Y_{v,t}(x)dv,\int_0^t b''(Y_{v,t}(x))D Y_{v,t}(x)dv \Big\rangle_{\mathcal{H}} \\
&& \geq e^{-2T\|b'\|_{\infty}}\int_0^T\Big(\int_0^t b''(v,Y_{v,t}(x)) \textbf{1}_{[v,t]}(\alpha) \ e^{- \int_{v}^{\alpha} b'(r, Y_{r,t}(x)) dr}\ dv\Big)^2 \ d\alpha \\
&& \geq  C^2 e^{-4T\|b'\|_{\infty}} \int_0^t \int_0^t\int_0^t \textbf{1}_{[v,t]}(\alpha) \textbf{1}_{[v_1,t]}(\alpha) dv_1 \ dv \ d\alpha \\
&& =  C^2 e^{-4T\|b'\|_{\infty}} \int_0^t \int_0^t (t-\max\{v_1,v\})  dv_1 \ dv \\
&& = C^2 e^{-4T\|b'\|_{\infty}} \int_0^t \int_0^t (t-\frac{v_1+v}{2}-\frac{|v_1-v|}{2}) dv_1 \ dv \\
&& = C^2 e^{-4T\|b'\|_{\infty}} \Big(\frac{t^3}{3}\Big) \\
&& >  0,\,\, \text{for all}\quad t>0\,,
\end{eqnarray*}
Thus, by the hypothesis \eqref{cc22}, for every $t\in (0,T]$, we see that
\begin{equation*}
A_2:=(u_{0} (Y_{t}(x)) ) ^{2}  \Vert DJY_{t}(x)\Vert_{\mathcal{H}} ^{2} \geq C^4 e^{-4T\|b'\|_{\infty}} \Big(\frac{t^3}{3}\Big):= C_5(t)>0\,.
\end{equation*}
Now, recalling that $\Vert Du(t,x) \Vert ^{2}_{\mathcal{H}} =A_{1}+A_{2}+A_{3}$ with $A_{1}>0$, $A_2\geq 0$ and $A_{3}\geq 0$ defined in the proof of Proposition \ref{existence}, we have, by the calculus above, that
\begin{equation}\label{6o-8}
\Vert Du(t,x) \Vert ^{2}_{\mathcal{H}}\geq C_5(t)>0\,.
\end{equation}
Therefore, the conclusion is obtained.\qed

\medskip


Using the Proposition \ref{Prob_B-C} (see \cite[Proposition 23]{Bally-Caramellino}), we can prove the H\"older continuity of the density.

\begin{prop}\label{regularity}  Assume  $b\in L ^{\infty} \left( (0,T), C^{3} _{b}(\mathbb{R})\right)$ such that
(\ref{cc11}) holds. Assume $u_{0} \in C_{b} ^{1} (\mathbb{R}) $  such that (\ref{cc22}) holds.

Let $u(t,x)$ be the solution to the stochastic continuity equation \eqref{continuity}. Then, for every $t\in [0,T]$ and for every $x\in \mathbb{R}$, the density $\rho_{u(t,x)}$ of the solution $u(t,x)$
  is H\"older continuous for any exponent $q<1$. Moreover, there exist  two universal constants $C,p>0$ such that
$$\rho_{u(t,x)}(z)\leq C \left( P (\vert u(t,x)-z\vert \leq 2) \right)^{\frac{1}{p}}$$
\end{prop}
{\bf Proof: } From Proposition \ref{p2} we have that $ u(t,x) $ belongs to $\mathbb{D} ^{2,p}$ for every $p\geq 2$. By Lemma \ref{l6}, we have that $\Vert Du(t,x) \Vert ^{-2}_{\mathcal{H}}$ belongs to $L ^{p}(\Omega)$ for every $t\in (0,T]$ and for every $x\in \mathbb{R}$ and for each $p\geq 1$. Then we can apply Proposition \ref{Prob_B-C}. \qed

\subsection{Gaussian estimates}

We finally prove the Gaussian estimate for the solution to the continuity equation.

\begin{prop}\label{p6}
Assume  $b\in L ^{\infty} \left( (0,T), C^{2} _{b}(\mathbb{R})\right)$ such that
(\ref{cc11}) holds. Assume $u_{0} \in C_{b} ^{1} (\mathbb{R}) $  such that (\ref{cc22}) holds. Let $(u(t,x))_{t\in [0,T], x\in \mathbb{R}}$ be the solution to the continuity equation \eqref{continuity}.
Then, for every $t\in  [t_{0},T]$ with $0 < t_{0} <T$  and for every $x\in \mathbb{R}$, the random  $u(t,x)$ admits a density $\rho _{u(t,x) }$ and  there exist two positive constants $c_{1}, c_{2}$ such that
\begin{equation}
\label{5a-5}
\frac{\mathbf{E} \vert u(t,x)- m\vert }{2c_{1} t} e ^{-\frac{(y-m) ^{2}}{2c_{2} t}}\leq \rho_{ u(t,x)} \leq  \frac{\mathbf{E} \vert u(t,x)- m\vert }{2c_{2} t}e ^{-\frac{(y-m) ^{2}}{2c_{1} t}}.
\end{equation}
\end{prop}
{\bf Proof: }
We will show that
\begin{equation}
\label{duu}
ct \leq \langle Du(t,x), \widetilde{Du(t,x)}\rangle _{\mathcal{H}}\leq Ct
\end{equation}
with $0<c<C$ where $\widetilde{Du(t,x)}$ is constructed from $Du(t,x)$ via the formula (\ref{tilde}). From (\ref{5o-5}),
\begin{eqnarray}
&&\langle Du(t,x), \widetilde{Du(t,x)}\rangle _{\mathcal{H}}\nonumber \\
&=& u_{0} '(Y_{t}(x)) \widetilde{u_{0} '(Y_{t}(x)) } JY_{t}(x) \widetilde{JY_{t}(x)} \langle DY_{t}(x), \widetilde{DY_{t}(x) } \rangle _{\mathcal{H}}
\nonumber \\
&& +  u_{0} (Y_{t}(x)) \widetilde{u_{0} (Y_{t}(x)) }
\langle DJY_{t}(x), \widetilde{DJY_{t}(x) } \rangle_{\mathcal{H}}\nonumber \\
&&+u_{0} '(Y_{t}(x))  JY_{t}(x)\widetilde{u_{0} (Y_{t}(x)) }\langle DY_{t}(x), \widetilde{DJY_{t}(x) } \rangle _{\mathcal{H}}\nonumber \\
&&+\widetilde{ u_{0} '(Y_{t}(x))} \widetilde{  JY_{t}(x)} u_{0} (Y_{t}(x)) \langle DJY_{t}(x), \widetilde{DY_{t}(x) } \rangle _{\mathcal{H}} \nonumber\\
&:=& I_{1}+ I_{2}+ I_{3}+ I_{4}.   \label{5o-10}
\end{eqnarray}
Note that, since $DY_{t}(x), \widetilde{DY_{t}(x)}, DJY_{t}(x) $ and $\widetilde{DJY_{t}(x)} $ are all negative, we can replace $\langle Du(t,x), \widetilde{Du(t,x)}\rangle _{\mathcal{H}}$ by $\langle \vert Du(t,x)\vert, \vert \widetilde{Du(t,x)}\vert \rangle _{\mathcal{H}}$  in (\ref{5o-10}) (and similarly for the other inner products in (\ref{5o-10})). 
We also have, from (\ref{5a-4}), for $t, \alpha \in (0, T], x\in \mathbb{R}$,
\begin{equation}\label{5o-7}
\vert D_{\alpha} Y _{t} (x)\vert =\vert 1_{[0,t]} (\alpha )e^{- \int_{0} ^{\alpha} b'(r, Y_{r,t}) dr}\vert \geq C1_{[0,t]} (\alpha ) e^{-\Vert b'\Vert _{\infty} \alpha}
\end{equation}
and  by (\ref{jy}),
\begin{equation}
\label{5o-8}
\vert JY_{t}(x) \vert  \geq  e ^{-\Vert b'\Vert _{\infty} t}
\end{equation}
while from (\ref{5o-6}),
\begin{equation}
\label{5o-9}
\vert D_{\alpha} JY_{t}(x)\vert \geq C t e^{\Vert b'\Vert_{\infty}t}1_{[0,t]} (\alpha ) e^{-\Vert b'\Vert_{\infty} \alpha}.
\end{equation}

Under  the hypothesis $b\in L ^{\infty} \left( (0,T), C^{2} _{b}(\mathbb{R})\right)$ and
(\ref{cc11}),  the  lower inequalities (\ref{5o-7}), (\ref{5o-8}) and (\ref{5o-9})  hold also for $\widetilde{DY_{t}(x)}$, $ \widetilde{JY_{t}(x)}$ and $ \widetilde{DJY_{t}(x)}$.  Consequently, since $t\in [t_{0}, T] $ with $0<t_{0}<T$ the summand

$$ \langle DJY_{t}(x), \widetilde{DJY_{t}(x) } \rangle_{\mathcal{H}}$$ is  bigger than $cT$ with $c>0$ and so $ I_{2} \geq ct$, via (\ref{cc22}). Since $I_{1}, I_{3}, I_{4} $ are positive we obtained the lower bound in (\ref{duu}).  By Lemmas \ref{l1}, \ref{l2} and \ref{l3}, the upper bound in (\ref{duu}) clearly holds. Therefore the conclusion is obtained.  \qed

%
%
%

{\bf Acknowledgement: } C. Olivera and C. Tudor  acknowledge partial support from the CNRS-FAPESP grant 267378.
C.  Olivera is partially supported by FAPESP   by the grants 2017/17670-0 and 2015/07278-0.


\end{document}